\documentclass[12pt]{article}
\usepackage{amssymb,amsmath,amscd,amsthm}
\usepackage[cp1251]{inputenc}
\usepackage[T2A]{fontenc}
\usepackage{pstricks}
\usepackage{graphicx,psfrag}

\usepackage[english]{babel}

\usepackage{cases}
\renewcommand{\setminus}{\diagdown}
\renewcommand{\emptyset}{\varnothing}

\newtheorem{theorem}{Theorem}[section]

\newtheorem{lemma}[theorem]{Lemma}

\date{}

\begin{document}

\title{Weyl type bound on positive Interior Transmission Eigenvalues }

\author{ E.Lakshtanov\thanks{Department of Mathematics, Aveiro University, Aveiro 3810, Portugal.   This work was supported by {\it FEDER} funds through {\it COMPETE}--Operational Programme Factors of Competitiveness (``Programa Operacional Factores de Competitividade'') and by Portuguese funds through the {\it Center for Research and Development in Mathematics and Applications} (University of Aveiro) and the Portuguese Foundation for Science and Technology (``FCT--Fund\c{c}\~{a}o para a Ci\^{e}ncia e a Tecnologia''), within project PEst-C/MAT/UI4106/2011 with COMPETE number FCOMP-01-0124-FEDER-022690, and by the FCT research project
PTDC/MAT/113470/2009 (lakshtanov@rambler.ru).} \and
B.Vainberg\thanks{Department
of Mathematics and Statistics, University of North Carolina,
Charlotte, NC 28223, USA. The work was partially supported   by the NSF grant DMS-1008132 (brvainbe@uncc.edu).}}

\maketitle

\begin{abstract}
This paper contains a lower bound of the Weyl type on the counting function of the positive eigenvalues of the interior transmission eigenvalue problem which justifies the existence of an infinite set of positive interior transmission eigenvalues.
We consider the classical transmission problem as well as the case where the inhomogeneous medium contains an obstacle. One of the essential components of the proof is an estimate for the D-t-N operator for the Helmholtz equation for positive $\lambda$ that replaces the standard parameter-elliptic estimate valid outside of the positive semi-axis.
\end{abstract}

\textbf{Key words:}
interior transmission eigenvalues, counting function, Weyl formula


\section{Introduction}
Interior transmission eigenvalues (ITEs) were introduced in the middle of 1980s  and soon became a classical object in the scattering theory, see, e.g., the recent review \cite{HadCak}. Their importance is based on the relation of ITEs to the far-field operator: if real $\lambda=k^2$ is not an ITE, then the far-field operator with the wave number $k$
 is injective and has a dense range. In particular,  when the linear sampling method (widely used in the inverse scattering theory) is applied for recovery of the support of the inhomogeneity in the medium, one needs to know  that the far-field operator has a dense range, i.e., $\lambda=k^2$ is not an ITE.  For this and other applications, it is important to know not only the discreteness of the ITEs but also their distribution. Note that ITEs can be measured, and this opens an opportunity to use ITEs for recovering properties of the scatterer (eg \cite[Th.3.2]{HadCak}).

Let us recall the definition of an ITE.
A value of $\lambda\in\mathbb C$ for which the homogeneous problem
\begin{numcases}{}
-\Delta u-\lambda u=0, & $x\in\mathcal O, \quad u\in H^2(\mathcal O)$\label{Anone0}\\
&\nonumber\\
-\nabla a(x)\nabla v-\lambda n(x)v=0, & $x\in\mathcal O, \quad v\in H^2(\mathcal O)$\label{Anone}\\
&\nonumber\\
\begin{array}{l}u-v=0, \\ \frac{\partial u}{\partial\nu}-a(x)\frac{\partial v}{\partial\nu}=0,\end{array} & $x\in\partial\mathcal O$\label{Antwo}
\end{numcases}
has a non-trivial solution is called an \emph{interior transmission eigenvalue}. Here $a(x),n(x),$ $x\in\overline{\mathcal O}$, are  smooth positive functions, $\nu$ is the outward unit normal vector, and $\mathcal O\subset\mathbb R^d$ is a bounded domain with a $C^\infty$-boundary.

The problem \eqref{Anone0}-\eqref{Antwo} appears naturally when the scattering of plane waves is considered and the inhomogeneity in $\mathbb R^d$ is located in $\mathcal O$ and described by functions $a(x),n(x)$. Our discussion primarily covers the cases $d=2,3$, but every result below can be automatically carried over to higher dimensions.

We will also consider the case where $\mathcal O$ contains a compact obstacle $\mathcal V \subset \mathcal O$, $\partial \mathcal V \in C^\infty$.
In this case, equation \eqref{Anone} is replaced by
\begin{equation}\label{AnoneB}
-\nabla a(x)\nabla v-\lambda n(x)v=0, \quad
x\in\mathcal O\setminus\mathcal V, \quad
v\in H^2(\mathcal O\setminus\mathcal V);~~~v(x)=0, \quad
x\in\partial\mathcal V,
\end{equation}
while equation \eqref{Anone0} remains valid in $\mathcal O$. For simplicity of notation, we will consider problem \eqref{Anone0}-\eqref{Antwo} as a particular case of \eqref{Anone0},\eqref{AnoneB},\eqref{Antwo} with $\mathcal V=\emptyset$. In \eqref{AnoneB} --- and in our previous papers on ITEs --- the Dirichlet boundary condition on $\partial\mathcal V$ can be replaced by the Neumann or Robin boundary condition without any changes in the results or proofs.

There are weaker definitions of ITEs where only $u-v\in H^2(\mathcal O\setminus\mathcal V)$ while $u\in L^2(\mathcal O), v\in L^2(\mathcal O\setminus\mathcal V)$, or all three functions $u,v$ and
 $u-v$ are assumed to be square integrable (the boundary conditions \eqref{Antwo} would still be meaningful so long as $u$ and $v$ satisfy the homogeneous elliptic equations). In fact, these weak eigenfunctions of the ITE problem belong to the Sobolev space $H^2$ under conditions imposed in the present paper (see the Attachment in \cite{lakvain7}), i.e., the a priory assumption $u,v\in H^2$ does not reduce the set of ITEs.

Denote the set of real {\it non-negative} ITEs with their multiplicities taken into account by $\{\lambda^T_i\}$. Similarly, denote the set of positive eigenvalues of the Dirichlet problem for $-\Delta$ in $\mathcal O$ by $\{\lambda_i\}$, and the set of positive $\lambda>0$ for which equation \eqref{AnoneB} in $\mathcal O\setminus V$ with the Dirichlet boundary conditions at the boundary $\partial(\mathcal O\setminus V)$ has a nontrivial solution by $\{\lambda^{a,n}_i\}$.
The corresponding counting functions will be denoted by
\begin{equation}\label{countdef}
N_T(\lambda)=\#\{i: \alpha < \lambda^T_i\leq\lambda\},
\quad
N(\lambda)=\#\{i:\lambda_i\leq\lambda\},
\quad
N_{a,n}(\lambda)=\#\{i:\lambda^{a,n}_i\leq\lambda\},
\end{equation}
where $\alpha $ is an arbitrary number from the interval $(0, \min(\lambda_0,\lambda^{a,n}_0))$ that does not belong to the set $\{ \lambda_i^T\}$. Note that we do not count positive ITEs on the segment $[0,\alpha]$ where we can not always guarantee that the number of ITEs is finite if $\mathcal V \neq \emptyset$.

Let us stress that problem \eqref{Anone0},\eqref{AnoneB},\eqref{Antwo} is not symmetric, and the existence of the real eigenvalues can not be established by soft arguments. Note, that even if there are only countably many positive ITEs, they could be distributed so sparsely or so densely, that from a practical point of view, the situation would be the same as whenever the set of ITEs is finite or non-discrete, respectively. Thus it is important to know conditions for the set  $\{\lambda^T_i\}$ to be discrete (counter examples can be found in \cite{lakvain5}), to be infinite, as well as to know the asymptotic behavior of $N_T(\lambda)$ as $\lambda\to\infty$.

The main result of the paper is stated in the following theorem.
\begin{theorem}\label{th} Assume that the set of ITEs is discrete, i.e., each ITE is isolated, (see Lemmas \ref{lemma1},\ref{lemma2} for sufficient conditions) and either $a(x)-1\neq0$, $x\in\partial\mathcal O$, or $a(x)\equiv1$ on $\mathcal O\setminus\mathcal V$, $n(x)-1\neq0$, $x\in\partial\mathcal O$. Let
\begin{equation*}
\gamma:=\operatorname{Vol}(\mathcal O)-\int_{\mathcal O\setminus\mathcal V}\left(\frac{n(x)}{a(x)}\right)^{\frac{d}{2}}dx\neq0.
\end{equation*}
Then the set of positive ITEs is infinite, and moreover,
\begin{equation}\label{theorem}
N_T(\lambda)   \geq \frac{\omega_d}{(2\pi)^d} |\gamma|   \lambda^{\frac{d}{2}} + O(\lambda^{\frac{d}{2}-\delta}), \quad \lambda \rightarrow \infty,
\end{equation}
where $\delta=\frac{1}{d+1}$ if $a(x) -1 \neq 0, x \in \partial \mathcal O$ and $\delta=\frac{1}{2d}$ if $a(x) \equiv 1,$  $\omega_d$ is the volume of the unit ball in $\mathbb R^d$.
\end{theorem}
When $\mathcal V=\emptyset$, the validity of the estimate $N_T(\lambda)   \geq C   \lambda^{d/2},~\lambda \rightarrow \infty,$ for some constant $C$ was obtained in \cite{ss} in the case where $a(x)\equiv1$ and $n(x)>1$ everywhere inside $\mathcal O$.

Theorem \ref{th} is a completion of our previous results \cite{lakvain6},\cite{lakvain7} where inequality \eqref{theorem} was proven for $\gamma$ having a specific sign, i.e., $\sigma\gamma>0$ where
\begin{equation*}
\sigma=\left\{\!\!\!
\begin{array}{ll}
\operatorname{sgn}(1-a(x))
&\text{if $a(x)\neq1$ for all $x\in\partial\mathcal O$,}\\
\operatorname{sgn}(n(x)-1)
&\text{if $a(x)\equiv1$ on $\mathcal O\setminus\mathcal V$ and $n(x)\neq1$ for $x\in\partial\mathcal O$.}
\end{array}\right.
\end{equation*}
\textit{Some previous results when $a(x)\neq1$ for $x\in\partial\mathcal O$.} Theorem \ref{th} requires the set of ITEs to be discrete.  The next lemma contains some known sufficient conditions which guarantee this property.
We call a set of eigenvalues discrete if each of them is isolated and has a finite multiplicity.
\begin{lemma}\label{lemma1}\textnormal{\cite{haddar1},\cite{LV4}}
\begin{enumerate}
\item If $(a(x)-1)(a(x)n(x)-1) \neq 0,~ x \in \partial \mathcal O$ then the set of ITEs is discrete.
\item If $a(x) \! - \!1 \neq 0, ~ x \in \partial \mathcal O$ and in \eqref{AnoneB}, $n(x)$ is replaced by $cn(x)$, then, for all
but finitely many $c\in \mathbb R$, the set of ITEs is discrete.
\item If $a(x)\!- \!1 \neq 0, ~ x \in \overline{\mathcal O},~\mathcal V=\emptyset$ and $\int_{\mathcal O} n(x) dx \neq 1$ then the set of ITEs is discrete.
\end{enumerate}
\end{lemma}
It was known (due to F.Cakoni, D.Gintides, H.Haddar, A.Kirsch, e.g., \cite[th.4.8]{HadCak},\cite{Dr},\cite{CK}) that the set $\{\lambda^T_i\}$ of non-negative ITEs is infinite provided that $\mathcal V = \emptyset$, and $(1-n(x))(a(x)-1)>0$  for each $x\in \overline{\mathcal O}$.

For $\mathcal V\neq\emptyset$ with $\operatorname{Vol}(\mathcal V)$ sufficiently small, the existence of at least one real eigenvalue whenever $a(x)>1$ and $n(x)\neq1$ on  $\overline{\mathcal O}$ was shown in \cite{chobst}. The authors of \cite{chobst} noted that the case $\mathcal V\neq\emptyset$ with $a(x)<1$ remains unstudied.

The asymptotics for the counting function of complex ITEs $\widehat{\lambda}^T_i$ is obtained in \cite{lakvain5} for the case $\mathcal V = \emptyset$ under condition 1 of Lemma \ref{lemma1}:
$$
\#\{i:|\widehat{\lambda}^T_i| \leq \lambda\} = \lambda^{\frac d 2} \frac{\omega_d}{(2\pi)^d} \int_{\mathcal O} \left (1+\left ( \frac{n(x)}{a(x)} \right )^{\frac d 2} \right ) dx + o(\lambda^{\frac d 2}), \quad \lambda \rightarrow \infty.
$$
Thus the right-hand side of this equality provides an  upper estimate on $N_T(\lambda) $ for real $\lambda^T_i$. It was also shown in \cite{lakvain5} that $\lim_{i \rightarrow \infty} 2\arg (\lambda^T_i)=0$.

\textit{Some previous results in the case $a(x)\equiv 1, x\in\mathcal O$ and $n(x)\neq1, x\in\partial\mathcal O$.} A sufficient condition for the discreteness of ITEs is given by the following lemma.
\begin{lemma}\label{lemma2}\textnormal{\cite{sylv},\cite{lakvain7},\cite{Rob},\cite{chobst}}
\begin{enumerate}
\item If $\mathcal V = \emptyset$ and $n(x)-1 \neq 0, \quad x \in \partial \mathcal O,$ then the set of  non-zero ITEs is discrete (i.e., each of them is isolated and has a finite multiplicity.) Point $\lambda=0$ is an isolated ITE of infinite multiplicity.
\item If $\mathcal V \neq \emptyset$  and  $n(x)-1 \neq 0, \quad x \in \partial \mathcal O,$ then the set of ITEs is at most countable with the only possible accumulation points at zero and infinity. All the ITEs have finite multiplicity and $\lambda=0$ is not an eigenvalue. 
\item If $\mathcal V \neq \emptyset$ and  $n(x) <1 , \quad x \in \overline{\mathcal O},$ then the set of non-zero ITEs is discrete.
\end{enumerate}
\end{lemma}
In \cite{Dr} it was proven that the set of the positive ITEs is discrete and infinite whenever $\mathcal V = \emptyset$ and the function $n(x)\!-\!1$ is nonzero for all $x\in \overline{\mathcal O}$.  More general conditions implying discreteness were obtained in \cite{ss}.

The existence of infinitely many real ITEs was proven in \cite{chobst} if $n(x)<1$ on $\mathcal O\setminus\mathcal V$.

It was shown in \cite{tsm} that whenever $\mathcal V=\emptyset$ and $n(x)>1$ on $\overline{\mathcal O}$, all but finitely many complex transmission eigenvalues are confined to a parabolic neighborhood of the real positive semi-axis. In \cite{tsm2} the authors justified the completeness of the set of the interior transmission eigenfunctions under the same assumption on $n(x)$. Later in \cite{Rob} the completeness of the eigenfunctions and the upper estimate by $C|\lambda|^{2+d/2}$ for the set of all ITEs was justified under the condition $\mathcal V=\emptyset$ with $n(x)\neq1$ on $\partial\mathcal O$.

\section{Proof of Theorem \ref{th}}

We use the approach developed in \cite{lakvain6},\cite{lakvain7} to prove the statement of the theorem in the case where $\sigma\gamma>0.$ In the first two steps below we outline the main parts of that proof which we need for studying the opposite sign of $\gamma$, and we refer to the original papers for more details.

\textit{Step 1. Reduction to the D-to-N operators.} We replace the ITE problem by the equivalent problem of finding values of $\lambda\geq 0$ for which the operator
\begin{equation}\label{opdn}
P(\lambda):=\sigma\big(F(\lambda)-F_{a,n}(\lambda)\big):H^{3/2}(\partial\mathcal{O})\to H^{3/2-s}(\partial\mathcal{O})
\end{equation}
has non-trivial kernels --- understood in a broader sense than usual (see Lemma \ref{ites} below). Here $H^{m}(\partial\mathcal{O})$ is the Sobolev space, $s$ is the order of the operator \eqref{opdn} which will be specified later, $F(\lambda)$ is the Dirichlet-to-Neumann map for equation \eqref{Anone0}, $F_{a,n}(\lambda)$ is the Dirichlet-to-Neumann map for problem \eqref{AnoneB} which is defined as follows:
\begin{equation}\label{DN}
F_{a,n}(\lambda):\phi\to\Big.\frac{\partial v}{\partial\nu}\Big|_{\partial\mathcal O},
\end{equation}
where $v$ satisfies \eqref{AnoneB} and
\begin{equation}\label{vvv}
v(x)=\phi, \quad x \in \partial \mathcal O.
\end{equation}
Note that $v$ vanishes at $\partial\mathcal V$, and functions from the domain of $F_{a,n}(\lambda)$ are defined only on $\partial \mathcal O$, not on $\partial (\mathcal O\setminus\mathcal V)$. The operators $F(\lambda),F_{a,n}(\lambda)$ are  symmetric elliptic pseudodifferential operators of the first order which are meromorphic in $\lambda\in\mathbb R$ with poles of the first order at the eigenvalues  $\lambda_i,\lambda_i^{a,n}$ of the Dirichlet problem for equations \eqref{Anone0},\eqref{AnoneB}, respectively.  If $a(x)\neq1$ on $\partial\mathcal O$, then the difference \eqref{opdn} remains to be a symmetric elliptic pseudodifferential operator of the first order (with a positive principal symbol) and, in this case, $s=1$ (see \cite{lakvain6}). However, the main symbol of the difference $F-F_{a,n}$ vanishes whenever $a(x)\equiv1$ on $\mathcal O\setminus\mathcal V$. If additionally, $n(x)\neq1$ on $\partial\mathcal O$, then $P(\lambda)$ is a symmetric elliptic pseudodifferential operator of order $-1$ with a positive principal symbol (see Lemma 1.1 of \cite{lakvain7} ), and $s=-1$.

The relationship between ITEs and the operator \eqref{opdn} is given by the next lemma which is a direct consequence of the definition of ITEs. Note that the phrase \emph{kernel of the operator} \eqref{opdn} is used below not only when the operator is analytic at $\lambda=\lambda_0$, but also when it has a pole at $\lambda=\lambda_0$. In the latter case, the kernel consists of the set of functions which are mapped to zero by both the residue and the principal part of the operator.
\begin{lemma}\label{ites}\textnormal{\cite{lakvain6},\cite{lakvain7}}
A point $\lambda=\lambda_0$ is an ITE if and only if the operator $F(\lambda_0) - F_{a,n}(\lambda_0)$ has a non trivial kernel or the following two conditions hold:
%
%
%
\begin{enumerate}
\item $\lambda=\lambda_0$ is an eigenvalue of the Dirichlet problem for $-\Delta$ and for equation \eqref{AnoneB}, i.e., $\lambda=\lambda_0$ is a pole for both $F(\lambda)$ and $F_{a,n}(\lambda)$.
\item The ranges of the residues of operators $F(\lambda)$ and $F_{a,n}(\lambda)$ at the pole $\lambda=\lambda_0$ have a non trivial intersection.
\end{enumerate}
Moreover, in each case, the multiplicity of the interior transmission eigenvalue $\lambda=\lambda_0$ is equal to $m_1+m_2$ where $m_1$ is the dimension of the kernel of the operator $F(\lambda)-F_{a,n}(\lambda)$ and $m_2$ is the dimension of the intersection of the ranges of the residues  at the pole $\lambda=\lambda_0$ ($m_2=0$ if $\lambda=\lambda_0$ is not a pole).
\end{lemma}

If $\lambda=\lambda_0$ satisfies the latter two conditions, we will call it a singular spectral point. Thus, singular spectral points belong to the intersection of three spectral sets: $ \{\lambda^{T}_i\},~\{\lambda_i\}$ and $\{ \lambda^{a,n}_i\}$.

\textit{Step 2. Relation between $N_T(\lambda)$
and the negative spectrum of the operator $P(\lambda)$.} The following approach is used to count the dimensions of the kernels of the operator \eqref{opdn} when $\lambda\to\infty$. It is based on studying the eigenvalues $\{\mu_i(\lambda)\}$ of the operator \eqref{opdn}.  Under the conditions imposed in Theorem \ref{th}, the operator $P(\lambda)$, $\lambda>0$, is an elliptic pseudodifferential operator of order $s$ with a positive principal symbol. For each $\lambda>0$, it has at most a finite number of negative eigenvalues. We denote this number by $n^-(\lambda) \geq 0$. Eigenvalues $\mu_i(\lambda)$ of the operator $P(\lambda)$ are meromorphic in
$\lambda \in \mathbb R$ with the possible poles only at points where $P(\lambda)$ has a pole.

Let us evaluate the difference $\widehat{n}(\lambda'):=n^-(\lambda')-n^-(\alpha)$ by moving $\lambda$ from $\lambda=\alpha$ (where $\widehat{n}=0$) to a fixed value $\lambda=\lambda'>\alpha$ which is not a pole of $P(\lambda)$. Here $\alpha>0$ is the constant defined in (\ref{countdef}). The eigenvalues $\mu_i(\lambda)$ may enter/exit the negative semi-axis $(-\infty,0)$ only through the end points of the semi-axis. Thus we can split $\widehat{n}(\lambda')$ as $\widehat{n}(\lambda')=n_1(\lambda')+n_2(\lambda')$ where $n_1(\lambda')$ is the number of eigenvalues $\mu_i(\lambda)<0$ that enter/exit the negative semi-axis $(-\infty,0)$ through the point $\mu=-\infty$ (when $\lambda$ changes from $\lambda=\alpha$ to $\lambda'>\alpha$) and $n_2(\lambda')$ is the number of eigenvalues $\mu_i(\lambda)<0$ that enter/exit the negative semi-axis $(-\infty,0)$ through the point $\mu=0$.
For example, $n_2(\lambda')>0$ if more eigenvalues enter the semi-axis through the origin than exit this semi-axis through the same point.

If an eigenvalue $\mu_i(\lambda)$ passes through the point $\mu=0$ when $\lambda$ passes through the value $\lambda_0\in(\alpha,\lambda')$, then $\lambda_0$ is an ITE due to Lemma \ref{ites}. Thus $N_T(\lambda)\geq |n_2(\lambda)|$. In fact, $n_2(\lambda)$ does not take into account singular ITEs, and we have
\begin{equation}\label{NT}
N_T(\lambda)\geq |n_2(\lambda)|+R(\lambda), \quad \lambda > \alpha,
\end{equation}
where $R(\lambda)$ is the counting function for the singular ITEs  located in the interval $(\alpha,\lambda)$. There are a couple of reasons why we can not claim equality in \eqref{NT}. In particular, $n_2(\lambda)$ does not count cases where $\mu_i(\lambda)$ vanishes at $\lambda=\lambda_0$ without crossing the origin $\mu=0$. A passage of an eigenvalue $\mu_i(\lambda)$ through $\mu=-\infty$ occurs when $\lambda$ pases through the poles of $P(\lambda)$. The poles are related to the eigenvalues of the Dirichlet problem for equations  \eqref{Anone0} and \eqref{AnoneB}. It was shown (in equation (14) of \cite{lakvain7}) that
\[
\big|n_1(\lambda)-\sigma\big(N_{a,n}(\lambda)-N(\lambda)\big)\big|\leq R(\lambda), \quad \lambda > \alpha,
\]
where $N(\lambda),N_{a,n}(\lambda)$ are counting functions for the eigenvalues of the Dirichlet problem for equations  \eqref{Anone0} and \eqref{AnoneB}, respectively. Thus
\begin{equation}\label{nnn}
n^-(\lambda)-n^-(\alpha) \geq \sigma\big(N_{a,n}(\lambda) - N(\lambda)\big) + n_2(\lambda)-R(\lambda), \quad \lambda > \alpha,
\end{equation}
and
\begin {equation}\label{nnnB}
n^-(\lambda)-n^-(\alpha) \leq \sigma\big(N_{a,n}(\lambda)-N(\lambda)\big) + n_2(\lambda)+R(\lambda), \quad \lambda > \alpha.
\end{equation}
These inequalities can be combined with the standard Weyl asymptotics for the counting functions of the Dirichlet problems:
\[
\sigma\big(N_{a,n}(\lambda) -N(\lambda) \big) = \frac{-\omega_d}{(2\pi)^d} \sigma \gamma   \lambda^{d/2}+O( \lambda^{(d-1)/2}), \quad \lambda\to\infty.
\]

Only inequality \eqref{nnnB}, but not \eqref{nnn} was used in \cite{lakvain6}, \cite{lakvain7}.
There it was assumed that  $\sigma\gamma>0$, i.e., $\sigma(N_{a,n}-N)\sim -\sigma\gamma\lambda^{d/2}<<0$. In this case, since $ n^-(\lambda)\geq0$, the statement of Theorem \ref{th} follows immediately from \eqref{NT} and \eqref{nnnB}.

Assume now that $\sigma \gamma<0$, i.e., $\sigma (N_{a,n}-N)\sim -\sigma \gamma \lambda^{d/2} >>0$. Then \eqref{NT} and  \eqref{nnn} imply
\begin{equation}\label{nnn4}
\begin{split}
N_T(\lambda)\geq-n_2(\lambda)+R(\lambda)
&\geq\sigma\big(N_{a,n}(\lambda) - N(\lambda)\big)-\big(n^-(\lambda)-n^-(\alpha)\big)\\
&=\frac{\omega_d}{(2\pi)^d}|\gamma|\lambda^{d/2}-\big(n^-(\lambda)-n^-(\alpha)\big)+O(\lambda^{(d-1)/2}).
\end{split}\end{equation}
The lack of an estimate on $n^-(\lambda)$ from above prevented us from proving the theorem earlier in the case where $\sigma\gamma <0$. We will now show that
\begin {equation}\label{nnn5}
 n^-(\lambda)= O(\lambda^{d/2-\delta}), \quad \lambda\to\infty,
\end{equation}
and this will complete the proof of Theorem \ref{th}.

\textit{Step 3. Proof of \eqref{nnn5}.}
We first consider the case where $a(x)\neq1$ on $\partial\mathcal O$. We write the operator \eqref{opdn} in the form
$$
P(\lambda)=\sigma\big[F(0)-F_{a,n}(0)\big]
+\sigma\big[\big(F(\lambda)-F(0)\big)+\big(F_{a,n}(\lambda)-F_{a,n}(0)\big)\big]
=:M+Q(\lambda)
$$
Consider a positive symmetric elliptic pseudodifferential operator $K$ on $\partial\mathcal O$ of order $1/2$. The operator $KP(\lambda)K=KMK+KQ(\lambda)K$ has the same number $n^-(\lambda)$ of negative eigenvalues as the operator $P(\lambda)$. We'll show that this number (for the operator $KP(\lambda)K$) satisfies \eqref{nnn5}.

According to the Weyl law  \cite[Th.1.6.1]{safvas} for the counting function of the eigenvalues of the Dirichlet problem, we have
\[
N_{a,n}(\lambda)=     \#\{i:\lambda^{a,n}_i \leq \lambda\}=c\lambda^{d/2}+O(\lambda^{\frac{d-1}{2}}), \quad \lambda\to\infty.
\]
From here it follows that the dimension $\kappa_1$ of the space $E_1=E_1(\lambda)$ spanned by the eigenfunctions whose eigenvalues lie within the interval $L:=(\lambda-\lambda^t,\lambda+\lambda^t)$, $1/2<t<1$, does not exceed $C\lambda^{d/2+t-1}$ as $\lambda\to\infty.$ One can easily specify the resolvent estimate in the proof of Lemma \ref{mmm} below (where operator $T(\lambda))$ is defined) and omit the eigenvalues $\lambda_i^{a,n}$ that belong to the interval $L$ from the right-hand side in \eqref{main1} if the operator $F_{a,n}$ is considered on the subspace $E_1^\bot$ orthogonal to $E_1$, i.e., the following analogue of \eqref{main1} is valid:
\begin{equation*}
\|T(\lambda)\phi\|_{H^{\frac{1}{2}}(\partial\mathcal O)} \leq C \lambda^{2-t}  \|\phi \|_{H^{-1/2}(\partial \mathcal O)}, \quad \phi\in E_1^\bot, \quad \lambda\to\infty.
\end{equation*}
A similar estimate (on a subspace $E_2^\bot$ of codimension $\kappa_2<C\lambda^{d/2+t-1}, ~\lambda\to\infty,$) is valid when $a(x)\equiv n(x)\equiv 1$. Thus there exists a constant $b>0$ such that
\begin{equation}\label{end}
\|KQ(\lambda)K\phi\|_{H^{0}(\partial\mathcal O)} \leq b \lambda^{2-t}  \|\phi \|_{H^{0}(\partial \mathcal O)}, \quad \phi\in E_1^\bot{\cap} E_2^\bot, \quad \lambda\to\infty.
\end{equation}

It was mentioned in step 1, that $M=P(0)$ is a symmetric elliptic operator of the first order on $\partial\mathcal O$ with a positive symbol. Thus $KMK$ is of order two, and from the Weyl law it follows that the counting function $N_{KMK}(\mu)$ for its eigenvalues has order $\mu^{\frac{m}{p}}$ where $m$ is the dimension of the manifold and $p$ is the order of the operator, i.e.,
\begin{equation}\label{wlaw}
N_{KMK}(\mu)=O(\mu^{\frac{d-1}{2}}), \quad \mu\to\infty.
\end{equation}
Thus the dimension $\kappa_3$ of the space $E_3=E_3(\lambda)\subset H^0(\partial\mathcal O)$ spanned by the eigenfunctions of $KMK$ with the eigenvalues smaller than or equal to $b \lambda^{2-t}$ does not exceed $C\left(\lambda^{2-t} \right )^{\frac{d-1}{2}}, ~\lambda\to\infty$. The estimate \eqref{end} implies that the following form is positive:
\begin{equation}\label{sqf}
(KP(\lambda)K\phi,\phi)>0, \quad \phi\in E_1^\bot{\cap} E_2^\bot{\cap} E_3^\bot\subset H^0(\partial\mathcal O).
\end{equation}
Therefore, $n^-(\lambda)\leq\kappa_1+\kappa_2+\kappa_3$. We choose $t$ in such a way that $\kappa_1=\kappa_2=\kappa_3$, i.e., $t=d/(d+1)$. This implies that $n^-(\lambda)=O(\lambda^{\frac{d}{2}-\frac{1}{d+1}})$, proving \eqref{nnn5} when $a(x)\neq1$ on $\partial\mathcal O$.

Now let us prove \eqref{nnn5} in the case where $a(x)\equiv1$ on $\mathcal O\setminus\mathcal V$ and $n(x)\neq1$ on $\partial\mathcal O$.
Formula \eqref{razl2} implies that
$$
P(\lambda) = \sigma\big(F(\lambda)-F_{a,n}(\lambda)\big)= D+\sigma \lambda[F'(0)-F_{a,n}' (0)] + \sigma[S_1(\lambda)-S(\lambda)] = :D+\lambda M+Q(\lambda),
$$
where $D=\sigma [F(0)-F_{a,n}(0)]$ and $S_1$ is the operator $S$ when $a(x)\equiv n(x)\equiv 1$. Obviously, $D=0$ if $\mathcal V=\emptyset$ and $D$ has an infinitely smooth integral kernel if $\mathcal V\neq\emptyset$. Consider the operator
$$
KP(\lambda)K = \lambda KMK+K[D+Q(\lambda)]K,
$$
where $K$ is an arbitrary positive elliptic pseudodifferential operator on $\partial\mathcal O$ of order $3/2$.

The operator $P(\lambda)$ was studied in \cite{lakvain7}.  It was proven there (see Lemma 1.1) that $P(\lambda)$, $0<|\lambda|\ll1,$ is an elliptic operator of order $-1$ with the principal symbol
\begin{equation*}\label{ddd}
\lambda\frac{|1-n(x)| }{2|\xi|}, \quad x \in \partial \mathcal O, \quad \xi \in \mathbb R^{d-1},
\end{equation*}
in any local coordinates such that the map from global to local coordinates is orthogonal at the boundary $\partial\Omega$. On the other hand, the operator $P(\lambda)$ is analytic at $\lambda=0$, and $D+Q(\lambda)$ has order $-3$. Thus $M$ is an elliptic operator of order $-1$ with a positive principal symbol (equal to $\frac{|1-n(x)|}{2|\xi| }$). Moreover, since the operator $P(\lambda)$ is symmetric, $M=P'(0)$ is also symmetric. Thus,
the Weyl law \eqref{wlaw} holds, and with an extra factor $\lambda$, we have
\begin{equation}\label{wlaw1}
N_{\lambda KMK}(\mu)=O\left ( \left ( \frac{\mu}{\lambda} \right )^{\frac{d-1}{2}}\right ), \quad \lambda>0, \quad  \mu\to\infty.
\end{equation}

From \eqref{main2} it follows that
\begin{equation*}
\|KQ(\lambda)K\|_{H^0} \leq C  \frac{\lambda^4}{\operatorname{dist}(\lambda,\{\lambda_i^{a,n}\})^2}, \quad \lambda \rightarrow \infty,
\end{equation*}
and for any $t, ~1/2<t<1,$ one can find (similarly to \eqref{end}) subspaces $E_1,E_2$ whose dimensions are equal to $\kappa=c\lambda^{d/2+t-1}$ such that
\begin{equation}\label{end2}
\|KQ(\lambda)K\phi\|_{H^{0}(\partial\mathcal O)} \leq b \lambda^{4-2t}  \|\phi \|_{H^{0}(\partial \mathcal O)}, \quad \phi\in E_1^\bot{\cap} E_2^\bot, \quad \lambda\to\infty.
\end{equation}
Since the operator $D$ is bounded, \eqref{end2} holds with $Q$ replaced by $Q_1=D+Q$.

From \eqref{wlaw1} it follows that the dimension $\kappa_3$ of the space $E_3=E_3(\lambda)\subset H^0(\partial\mathcal O)$ spanned by the eigenfunctions of $\lambda KMK$ whose eigenvalues are smaller than or equal to $b\lambda^{4-2t}$ does not exceed $C(\lambda^{3-2t})^{\frac{d-1}{2}}$ as $\lambda\to\infty$. Now choose $t$ in such a way that $\kappa=\kappa_3$, i.e., $t=1-\frac{1}{2d}$. Then the quadratic form \eqref{sqf} is positive on the subspace of codimension $C(\frac{d}{2}-\frac{1}{2d})$. This completes the proof of \eqref{nnn5}. \qed

\section{The main technical lemma}
Recall that $F_{a,n}$ is the Dirichlet-to-Neumann operator for problem \eqref{AnoneB} (defined in \eqref{DN}), and $\{\lambda_i^{a,n}\}$ is the set of eigenvalues for the Dirichlet problem for equation \eqref{AnoneB}. The following statement is the main technical lemma of the paper.
\begin{lemma}\label{mmm} 1) The following expansion is valid
$$
F_{a,n}(\lambda)= F_{a,n}(0)+T(\lambda),
$$
where $F_{a,n}(0)$ is a pseudodifferential operator on $\partial\mathcal{O}$ of the first order, operator
$T(\lambda)$ has order $-1$, and
\begin{equation}\label{main1}
\|T(\lambda)\phi\|_{H^{\frac{1}{2}}(\partial\mathcal O)} \leq C \left (\frac{\lambda^2}{\operatorname{dist}(\lambda,\{\lambda_i^{a,n}\})}+\lambda  \right )  \|\phi \|_{H^{-1/2}(\partial \mathcal O)}, \quad \lambda\geq 0.
\end{equation}

2) One can also write $F_{a,n}(\lambda)$ in the form
\begin{equation}\label{razl2}
F_{a,n}(\lambda)= F_{a,n}(0)+\lambda F_{a,n}'(0) + S(\lambda)
\end{equation}
where operator $F_{a,n}'(0)$ is of order $ -1$, operator $S(\lambda)$ has order $-3$ and
\begin{equation}\label{main2}
\|S(\lambda)\phi\|_{H^{\frac{5}{2}}(\partial\mathcal O)} \leq C \left (\frac{\lambda^2}{\operatorname{dist}(\lambda,\{\lambda_i^{a,n}\})}+\lambda  \right ) ^2 \|\phi \|_{H^{-1/2}(\partial \mathcal O)}, \quad \lambda\geq 0.
\end{equation}
\end{lemma}

\textbf{Proof.} It is well known that $F_{a,n}(\lambda),F_{a,n}(0)$ are pseudodifferential operators of the first order, see \cite{lakvain7} for more details and references. There, one can also find the calculation of the full symbol of the operator $F_{a,n}(\lambda)$ which shows that the operator $F_{a,n}(\lambda)-F_{a,n}(0)$ has order $-1$. This fact also follows from \eqref{main1}. Thus in order to prove the first statement of Lemma~\ref{mmm}, it is enough to justify \eqref{main1}.

Let $v$ be the solution of the problem \eqref{AnoneB},\eqref{vvv}, and let $w$ be the solution of the same problem with $\lambda=0$, i.e.,
\begin{equation}\label{bb}
[F_{a,n}(\lambda)-F_{a,n}(0)]\phi=\Big.\frac{\partial(v-w)}{\partial n}\Big|_{\partial\mathcal O}.
\end{equation}
We assume that $\phi\in H^{-1/2}(\partial \mathcal O)$ and $v,w\in H^0=H^0(\mathcal O\setminus\mathcal V)$. One can consider non-smooth solutions of homogeneous elliptic problems (see \cite{roit}). Non-smooth solutions $v$ are understood as limits in $H^0$ of solutions $v_n$ to the same problem with smooth boundary conditions $\phi_n=v_n|_{\partial \mathcal O}$ such  that $\phi_n\to\phi$ in  $H^{-1/2}(\partial \mathcal O)$ as $n\to\infty$. A solution $v$ exists and is unique if $\lambda$ is not an eigenvalue of the problem. Solutions $w$ are defined similarly. A~priori estimates are valid for these solutions (see \cite{roit}). In particular,
\begin{equation}\label{aprest}
\|w\|_{H^0} \leq C\|\phi\|_{H^{-1/2}(\partial\mathcal O)} .
\end{equation}

Since $v-w$ vanishes at the boundary and
\begin{equation}\label{ba}
-\nabla a(x) \nabla (v-w) - \lambda n(x)(v-w)=\lambda n(x)w, \quad x\in \mathcal O\setminus\mathcal V,
\end{equation}
the standard resolvent estimate implies that
\begin{equation}\label{ba1}
\|v-w\|_{H^0} \leq \frac{C\lambda}{\operatorname{dist}(\lambda,\{\lambda_i^{a,n}\})} \|w\|_{H^0}.
\end{equation}
Furthermore,
\begin{equation}\label{razl}
-\nabla a(x) \nabla (v-w)= \lambda n (v-w)+\lambda n w,  \quad x\in \mathcal O\setminus\mathcal V,
\end{equation}
and the operator $(-\nabla a(x) \nabla)^{-1}:H^0\to H^2$ is bounded. Thus
\begin{equation}\label{baba}
\|v-w\|_{H^2(\mathcal O \setminus\mathcal V)} \leq C \|\lambda(v-w)\|_{H^0}+\|\lambda w\|_{H^0} \leq C \left (\frac{\lambda^2}{\operatorname{dist}(\lambda,\{\lambda_i^{a,n}\})}+\lambda  \right ) \|w\|_{H^0},
\end{equation}
and therefore, due to \eqref{bb} and \eqref{aprest}, we have
$$
\left\|\frac{\partial}{\partial\nu}(v-w)\right\|_{H^{1/2}(\partial\mathcal O)}\leq C \left (\frac{\lambda^2}{\operatorname{dist}(\lambda,\{\lambda_i^{a,n}\})}+\lambda  \right )  \|\phi \|_{H^{-1/2}(\partial \mathcal O)},
$$
which justifies \eqref{main1}. This proves the first statement in Lemma~\ref{mmm}.

Let us now prove the second statement. First of all note that $F'_{n,a}(0)=T'(0)$ and therefore the operator $F'_{n,a}(0)$ has order $-1$. Thus, we only need to prove \eqref{main2}. Furthermore, it follows from \eqref{bb} that
\begin{equation*}
S(\lambda)\phi=\Big.\frac{\partial(v-w-\lambda z)}{\partial \nu}\Big|_{\partial\mathcal O}
\quad\text{where}\quad
z=\Big.\frac{\partial(v-w)}{\partial \lambda}\Big|_{\lambda=0}.
\end{equation*}
Differentiating \eqref{ba} in $\lambda$ and then putting $\lambda=0$, since $w$ does not depend on $\lambda$, we obtain
\begin{equation*}
-\nabla a(x) \nabla z - \lambda n(x)z= n(x)w+n(x)(v-w),\quad x\in \mathcal O\setminus\mathcal V.
\end{equation*}
Thus $u:=v-w-\lambda z$ satisfies
\begin{equation}\label{bab}
-\nabla a(x) \nabla u - \lambda n(x)u= -\lambda n(x)(v-w),\quad x\in \mathcal O\setminus\mathcal V;\quad u|_{\partial(\mathcal O\setminus\mathcal V)}=0.
\end{equation}
The same arguments that were used to get \eqref{baba} from \eqref{ba} allow us to obtain the following consequence of \eqref{bab}:
$$
\|u\|_{H^2(\mathcal O \setminus\mathcal V)}\leq C\left(\frac{\lambda^2}{\operatorname{dist}(\lambda,\{\lambda_i^{a,n}\})}+\lambda\right)\|v-w\|_{H^0}.
$$
We move the term $ \lambda n(x)u$ to the right-hand side of equation \eqref{bab} and use the a~priori estimate for the operator $\nabla a(x) \nabla$. Together with \eqref{ba1} and \eqref{baba}, the previous estimate leads to
$$
\|u\|_{H^4(\mathcal O\setminus\mathcal V)}\leq C\lambda\big(\|u\|_{H^2(\mathcal O\setminus\mathcal V)}+\|v-w\|_{H^2(\mathcal O\setminus\mathcal V)}\big)\leq C\left(\frac{\lambda^2}{\operatorname{dist}(\lambda,\{\lambda_i^{a,n}\})}+\lambda\right)^2\|w\|_{H^0},
$$
which together with \eqref{aprest} justifies \eqref{main2}. \qed

\textbf{Acknowledgment.}
The authors are grateful to L. Friedlander and Yu. Safarov for very useful discussions.

\end{document}